\documentclass[11pt, reqno]{amsart}
\theoremstyle{plain}
\newtheorem{theorem}{Theorem}[section]

\newtheorem{definition}[theorem]{Definition}

\begin{document}


\title[riccati inequality ]{Solvability of Matrix riccati inequalities  }
\thanks{Adviser: Professor Nikita Barabanov}
\author[K. Kissi]{Kevin Kissi}

\address{Department of Mathematics 2750\\ North Dakota State University\\PO BOX 6050\\ Fargo, ND 58108-6050\\ USA}

\email{kevin.kissi@ndsu.edu}

%
%

\begin{abstract}
 We consider matrix Riccati inequality arising in the theory of absolute stability, $H_\infty$ control problem, $LQ$ problem, and optimal estimation problem. In the case of sign definite frequency domain function, the solvability of Riccati inequalities is a subject of the famous Kalman- Yakubovich lemma. 
This paper presents necessary and sufficient conditions for solvability of Riccati inequality in the general sign indefinite case. To this end we use special representations of Hamiltonian matrices. The results are illustrated by an example. 
\end{abstract}

%
%
\maketitle

%
%

\section{Introduction and Problem statement}

Consider the Riccati inequality

\begin{equation}\label{1}
HA + A^{*}H + G - HB\Gamma ^{-1}B^{*}H < 0,
\end{equation}

\noindent where $A, B, G, \Gamma $ are given matrices of dimensions $n\times n, n\times m,
n\times n$, and $m\times m$ respectively; $G, \Gamma $ are Hermitian matrices, $\det \Gamma
\neq 0$. We are looking for necessary and sufficient conditions for existence of stabilizing and anti stabilizing solutions of the inequality (\ref{1}), that is, for matrices $H_{-}$ and $H_{+}$ such that (\ref{1}) holds for both of them, and matrices $A-B\Gamma^{-1}B^{*}H_{-}$, $-(A-B\Gamma^{-1}B^{*}H_{+})$ are Hurwitz.

In case when matrix $\Gamma$ is sign definite, the answer to this problem is given in the famous Kalman-Yakubovich lemma \cite{frequencydomain}. This lemma, has been originally proved in \cite{ineq}, then extended to infinite dimensional case in \cite{control_theory} and formulated in the most general form in \cite{solvability}.

For positive definite matrices $\Gamma$ solvability of inequality (\ref{1}) may be reduced to solvability of the algebraic Riccati equation (ARE)

$$ HA + A^{*}H + G - HB\Gamma ^{-1}B^{*}H = 0,$$

\noindent which is closely related to the existence and properties of maximal $J$-orthogonal invariant subspaces of Hamiltonian matrices (see \cite{ricci}, \cite{Hamiltonianmatrices}). Such equations are very important in the theories of optimal control, absolute stability, game theory, $H_{\infty}$ control (in the last two cases the matrix $\Gamma$ is sign indefinite) \cite{Scherer}. 

A number of important results concerning solvability of Riccati inequalities were published in \cite{Riccatiinequalities}. In this paper a new idea related to positive definiteness of matrix $M$ (see the text below) is used.

Theory of absolute stability study systems of the form 

\begin{equation}\label{abs} 
\begin{array}{c}
\frac{dx}{dt} = Ax + B\xi,\quad \sigma=C^{*}x, \\
 \xi = \varphi(\sigma,t)
\end{array}
\end{equation}

\noindent in the class $N_{F}$ of nonlinearities $\varphi$ satisfying the following Local Quadratic Constraint (LQC):

\begin{equation}\label{LQC} 
F(x,\varphi(C^{*}x,t))\ge 0 \quad \forall x,t,
\end{equation}

\noindent where $F$ is given quadratic form. The first and most known LQC is so-called sector condition:

$$ (\xi-\alpha C^{*}x)^{*}\Gamma (\beta C^{*}x-\xi)\ge 0,$$

\noindent where $\alpha$, $\beta$, $\Gamma>0$ are diagonal matrices, but there are a lot of other useful constraints \cite{frequencydomain}. 

System (\ref{abs}) is called absolutely stable in the class $N_{F}$, if system (\ref{abs}) with every function $\varphi\in N_{F}$ is globally asymptotically stable, and this stability is uniform with respect to functions $\varphi\in N_{F}$.

If there exists a positive definite quadratic form $V(x)=x^{*}Hx$ such that $V(x(t))$ is decreasing along all non zero solutions of systems (\ref{abs}) with all functions $\varphi\in N_{F}$, then system (\ref{abs}) is absolutely stable in $N_{F}$. 

Such Hermitian matrix $H$ exists if and only if the quadratic form

$$ x^{*}H(Ax+B\xi) <0 $$

\noindent for all non zero vectors $(x^{*},\xi^{*})^{*}$, for which $F(x,\xi)\ge 0$. 

According to Dines theorem such matrix $H$ exists if and only if there exists a Hermitian matrix $H$ such that

\begin{equation}\label{HF} 2x^{*}H(Ax+B\xi) + F(x,\xi) < 0
\end{equation}

\noindent for all non zero vectors $(x^{*},\xi^{*})^{*}$.

Assume $F(x,\xi)=x^{*}Gx-\xi^{*}\Gamma\xi$ with $\Gamma>0$. Then inequality (\ref{HF}) holds if and only if the following Riccati inequality holds

\begin{equation}\label{ric}
 HA + A^{*}H + G + HB\Gamma^{-1}B^{*}H < 0.
\end{equation}

Notice that the quadratic term of the left hand side is sign semidefinite.

Another example of problems where the Riccati inequalities arise is the problem of $H_{\infty}$ control.

Consider system

\begin{equation}\label{hinf}
\frac{dx}{dt} = Ax + B_{w}w + B_{u}u,
\end{equation}

\noindent where matrices $A$, $B_{w}$, $B_{u}$ are constant, $x$ is state vector, $w$ is exogenous input (noise for example), and $u$ is control. 

Consider a quadratic function 

$$ F(x,w,u) = x^{*}Gx + u^{*}\Gamma_{u}u - w^{*}\Gamma_{w}w,$$

\noindent where constant Hermitian matrices $G$, $\Gamma_{u}>0$ and $\Gamma_{w}>0$ are given. The problem consists of finding a controller $u=hx$ such that the closed loop matrix $A+B_{u}h$ is Hurwitz, and for all non zero functions $w\in L_{2}(0,\infty)$ along solutions $x$ with trivial initial value $x(0)=0$ we have

$$ \int_{0}^{\infty} F(x(t),w(t),u(t)) dt <0. $$

This problem may be reduced to solvability of the following 
Riccati inequality

$$ HA + A^{*}H + G + HB_{w}\Gamma_{w}^{-1}B_{w}^{*}H - HB_{u}\Gamma_{u}^{-1}B_{u}^{*}H < 0.$$

Notice that matrix of quadratic form in this inequality in many cases is sign indefinite.

Now consider a relation between the Riccati inequalities and the  Kalman-Yakubovich lemma.

Without loss of generality we assume that pair $(A,B)$ is controllable and matrix $A$ has no pure imaginary eigenvalues.

If matrix $\Gamma $ is negative definite, then inequality (\ref{1}) may be represented as
linear matrix inequality (LMI):

$$
\left(\begin{array}{cc}
HA + A^{*}H + G&HB \\
 B^{*}H&\Gamma \end{array}\right)  < 0. \; \;
$$

It may be solved via well-known technique.

Solvability of inequality (\ref{1}) is also a subject of the famous Kalman-Yakubovich 
lemma \cite{frequencydomain}. According to this lemma, inequality (\ref{1}) has a solution, 
if and only if, the following frequency domain inequality holds:

\begin{equation}\label{2}
\pi (i\omega ) < 0,
\end{equation}

\noindent for all $\omega \in [-\infty , \infty ]$, where

$$
\pi (\lambda ) = \Gamma  + B^{*}(\lambda I+A^{*})^{-1}G(A-\lambda I)^{-1}B,
$$

\noindent or, which is the same,

\begin{equation}\label{3}
\det \pi (i\omega ) \neq  0
\end{equation}

\noindent for all $\omega \in [-\infty, \infty ]$.

But if matrix $\Gamma $ is not sign definite, the inequality (\ref{3}) proved to be no
longer necessary for solvability of inequality (\ref{1}). In this paper we present the
desired necessary and sufficient conditions, which may be considered as generalization 
of Kalman-Yakubovich lemma to the case of sign-indefinite quadratic forms.

\section{Hamiltonian matrices}

In this section we present a new necessary condition for solvability of inequality
(\ref{1}), which will be proved later on to be also sufficient.

In the sequel we shall use the following matrices:

$$
R = \left(\begin{array}{cc}
A&-B\Gamma ^{-1}B^{*} \\
 -G&-A^{*}\end{array}\right) , \; \; \qquad J = \left(\begin{array}{cc}
0&-I \\
 I&0\end{array}\right).
$$

Matrix $JR$ is clearly Hermitian, therefore, matrix $R$ is ($J$-)Hamiltonian \cite{Hamiltonianmatrices}.

The set of eigenvalues of Hamiltonian matrix is symmetric with respect to the imaginary axis. Indeed, if $\det(R-\lambda I)=0$, then 

$$\begin{array}{c} 0 = \det(R-\lambda I)\det(J)=\det(JR-\lambda J) = \det(R^{*}J^{*}-\lambda J) = \\
=\det(R^{*}(-J)-\lambda J) = \det(-R^{*}-\lambda I)det(J),\end{array}$$

\noindent Therefore

$$ 0 = \overline{\det(R^{*}+\lambda I)} = \det(R-(-\bar\lambda I)),$$

\noindent and $-\bar\lambda$ is an eigenvalue of $R$.

\noindent It has been shown \cite{ricci} that the structure of Jordan blocks corresponding to eigenvalues $\lambda$ and $-\bar\lambda$ of matrix $R$ coincide. 

Assume matrix $R$ has no eigenvalues on the imaginary axis. Then there exist $n\times n$-matrices $\Lambda$, $X_{1}$, $\Psi_{1}$, $X_{2}$, $\Psi_{2}$ such that all eigenvalues of matrix $\Lambda$ have negative real parts, and 

\begin{equation}\label{RI}
 R \left(\begin{array}{cc} X_{1} & X_{2} \\ \Psi_{1} & \Psi_{2} \end{array}\right) =
\left(\begin{array}{cc} X_{1} & X_{2} \\ \Psi_{1} & \Psi_{2} \end{array}\right)
\left(\begin{array}{cc} \Lambda & 0 \\ 0 & -\Lambda^{*} \end{array}\right).
\end{equation}

\noindent Hence, we have the following equalities:

$$\begin{array}{c} 
AX_{1} - B\Gamma^{-1}B^{*}\Psi_{1} = X_{1}\Lambda \\ -GX_{1}-A^{*}\Psi_{1}=\Psi_{1}\Lambda. \end{array} $$

Assume matrix $X_{1}$ is nonsingular. Multiply the first equation by $\Psi_{1}X_{1}^{-1}$ from the left, by $X_{1}^{-1}$ from the right. Multiply the second equation by $X_{1}^{-1}$ from the right. Add the equations. Then with notation $H=\Psi_{1}X_{1}^{-1}$ we have

\begin{equation}\label{RE}
HA + A^{*}H + G - HB\Gamma^{-1}B^{*}H = 0.
\end{equation}

Hence, $H$ is a solution of the Riccati equation (\ref{RE}). 

Why $H$ is Hermitian? Denote $Z=col(X_{1},\Psi_{1})$. Then $RZ=Z\Lambda$ and

$$ Z^{*}JZ\Lambda =Z^{*}JRZ = Z^{*}(JR)^{*}Z = (RZ)^{*}J^{*}Z = (-\Lambda^{*})Z^{*}JZ.$$

\noindent But matrices $\Lambda$ and $-\Lambda^{*}$ have no common eigenvalues. Therefore 

$$ Z^{*}JZ=0. $$

\noindent Recalling the definition of $Z$ we get 

$$ -X_{1}^{*}\Psi_{2} + \Psi_{1}^{*}X_{1}=0,$$

\noindent which it turn is equivalent to $(\Psi_{1}X_{1}^{-1})^{*}=\Psi_{1}X_{1}^{-1}$, or $H^{*}=H$. 

Notice that $A-B\Gamma^{-1}B^{*}H=X_{1}\Lambda X_{1}^{-1}$. Therefore matrix $A-B\Gamma^{-1}B^{*}H$ is Hurwitz, and $H$ is a stabilizing solution of the Riccati equation (\ref{RE}). 

If we consider matrix $col(X_{2},\Psi_{2})$ instead of $Z$, we arrive to a solution $H$ such that matrix $-(A-B\Gamma^{-1}B^{*}H)$ is Hurwitz. Then $H$ is anti stabilizing solution of the Riccati equation (\ref{RE}).

Inverse, if $H$ is a solution of the Riccati inequality (\ref{RI}), then for some positive definite matrix $\Delta G$, matrix $H$ is a solution of the Riccati equation 

\begin{equation}\label{REM}
HA + A^{*}H + G + \Delta G - HB\Gamma^{-1}B^{*}H = 0.
\end{equation}

Assume $H$ is a stabilizing solution of equation (\ref{REM}). Then for some representation (\ref{RI}) of matrix 

$$ R_{new} = \left(\begin{array}{cc}
A & -B\Gamma^{-1}B^{*} \\  -G-\Delta G & -A^{*} \end{array}\right)
$$

\noindent we have $A-B\Gamma^{-1}B^{*}H=X_{1}\Lambda X_{1}^{-1}$, and therefore matrix $\Lambda$ is Hurwitz. In particular, it means that matrix $R_{new}$ has no pure imaginary eigenvalues.

Thus, the problem of finding necessary and sufficient conditions for existence of stabilizing and anti stabilizing solutions of Riccati inequality (\ref{1}) is reduced to a problem of existence of a positive definite matrix $\Delta G$ such that matrix $R_{new}$ has no pure imaginary eigenvalues, and in the representation (\ref{RI}) matrices $X_{1}$ and $X_{2}$ are nonsingular.

\section{Special Transformation}

Assume $V_{1}$, $V_{2}$ are $n\times m$-matrices, and $V=col(V_{1},V_{2})$. Denote by $R(V)$ the Hamiltonian matrix 

$$ R + V(JV)^{*} = \left(\begin{array}{cc} 
A-V_{1}V_{2}^{*} & -B\Gamma^{-1}B^{*}+V_{1}V_{1}^{*} \\  
-G-V_{2}V_{2}^{*} & -A^{*}+V_{2}V_{1}^{*} \end{array}\right).
$$

The corresponding Riccati equation has a form

\begin{equation}\label{REV}
H(A-V_{1}V_{2}^{*})+(A-V_{1}V_{2}^{*})^{*}H + G + V_{2}V_{2}^{*}-
HB\Gamma^{-1}B^{*}H + HV_{1}V_{1}^{*}H = 0.
\end{equation}

For a solution $H$ of this equation we have 

$$ HA + A^{*}H + G - HB\Gamma^{-1}B^{*}H = -(V_{2}-HV_{1})(V_{2}-HV_{1})^{*} \le 0.$$

\noindent If the right hand side is strictly negative, then $H$ is a solution of inequality (\ref{1}). 

\section{Special case}

Assume matrix $R$ has no pure imaginary eigenvalues. Then for sufficiently small positive number $\epsilon$ matrix $R(\epsilon I)$ also has no pure imaginary eigenvalues. 

Consider a representation (\ref{RI}) of matrix $R(\epsilon I)$. For every positive number $\delta$ there exists a Hamiltonian matrix 

$$ \tilde R = \left(\begin{array}{cc} 
\tilde A & -\tilde B\tilde\Gamma ^{-1}\tilde B^{*} + \epsilon I \\  
-\tilde G -\epsilon I & -\tilde A^{*} \end{array}\right)
$$

\noindent such that $\|R(\epsilon I)-\tilde R\|<\delta$, and in the representation (\ref{RI}) of matrix $\tilde R$ matrices $X_{1}$, $X_{2}$ are nonsingular. Denote by $H$ the stabilizing solution of corresponding Riccati equation. Then

$$\begin{array}{c} HA + A^{*}H + G - HB\Gamma^{-1}B^{*}H = \\
= H(A-\tilde A) + (A-\tilde A)^{*}H + (G-\tilde G+\epsilon I) - \\
- (HB\Gamma^{-1}B^{*}H - H\tilde B\tilde\Gamma^{-1}\tilde B^{*}H) + H\epsilon H.
\end{array}
$$

\noindent For sufficiently small number $\delta$ the right hand side of the equality is negative. Therefore, $H$ is a solution to the Riccati inequality (\ref{1}). Since matrix $\tilde A - \tilde B\tilde\Gamma^{-1}\tilde B^{*}H$ is Hurwitz, for sufficiently small number $\delta$ the matrix $A-B\Gamma^{-1}B^{*}H$ is Hurwitz, and $H$ is a stabilizing solution of inequality (\ref{1}).

The same conclusion is true for the anti stabilizing solution of inequality (\ref{1}).

Thus, if matrix $R$ has no pure imaginary eigenvalues, then the Riccati inequality (\ref{1}) has both stabilizing and anti stabilizing solutions.

\section{General case}

Now consider general case: matrix $R$ may have pure imaginary eigenvalues. Our next goal is to figure out how eigenvalues of matrix $R+V(JV)^{*}$ depend on $V$.

We need the following definitions. 

\begin{definition} Let $P$ be a subspace of ${\bf C}^{2n}$, and $T$ be a matrix such that set of columns of $P$ is a basis of $P$. Denote by $n_{-}(P)$, $n_{0}(P)$, and $n_{+}(P)$ the numbers of positive, zero, and negative eigenvalues of matrix $T^{*}iJT$ respectively. 
\end{definition}

Obviously, these numbers do not depend on the choice of basis $T$ of $P$.

Denote by $J_{1},\ldots,J_{m}$ all Jordan blocks of matrix $R$ with pure imaginary eigenvalues $i\omega_{1},\ldots,i\omega_{m}$ respectively:

$$ J_{j} = \left(\begin{array}{ccccc} i\omega_{j} & 1 & 0 & \ldots & 0 \\
0 & i\omega_{j} & 1 & \ldots & 0 \\ \ldots & \ldots & \ldots & \ldots & \ldots \\ 0 & 0 & 0 & \ldots & i\omega_{j} \end{array}\right).
$$

\noindent We assume that the Jordan blocks are arrange in such an order that $\omega_{1}\le\omega_{2}\le\ldots\le\omega_{m}$. Let $P_{1},\ldots,P_{m}$ be $R$-invariant subspaces associated to Jordan blocks $J_{1},\ldots,J_{m}$ of dimensions $n_{1},\ldots,n_{m}$ respectively. It is known that for every $j=1,\ldots,m$ the value $n_{+}(J_{j})-n_{-}(J_{j})$ is equal to $-1$, $0$, or $1$, and $n_{0}(J_{j})=0$. We use classification given by M. Krein.

\begin{definition} We say that the block $J_{j}$ contains $n_{+}(J_{j})$ eigenvalues $i\omega_{j}$ of the first type, and $n_{-}(J_{j})$ eigenvalues $i\omega_{j}$ of the second type.
\end{definition}

For every subspace $P_{j}$ there exist a number $\beta_{j}\in\{-1,1\}$, and a matrix $S_{j}$ such that the columns of $S_{j}$ span $P_{j}$, $RS_{j}=S_{j}J_{j}$, 

$$ S_{j}^{*}JS_{j} = \epsilon_{j} \left(\begin{array}{ccccc} 0 & 0 & \ldots & 0 & -1 \\
0 & 0 & \ldots & (-1)^{2} & 0 \\  \ldots &  \ldots &  \ldots &  \ldots &  \ldots \\ (-1)^{m_{j}} & 0 & \ldots & 0 & 0 \end{array}\right),
$$

\noindent and

(a) if the size $n_{j}$ of $J_{j}$ is even, then $n_{+}(J_{j})-n_{-}(J_{j})=0$ and $\epsilon_{j}=(-1)^{n_{j}/2}\beta_{j}$ ;

(b) if the size $n_{j}$ of $J_{j}$ is odd, then $n_{+}(J_{j})-n_{-}(J_{j})=\beta_{j}$ and $\epsilon_{j}=(-1)^{(n_{j}-1)/2}i\beta_{j}$.

\noindent Moreover, matrices $S_{j}$ for distinct $j$ are $J$-orthogonal: $S_{j_{1}}^{*}JS_{j_{2}}=0$ if $j_{1}\ne j_{2}$.

The value $\beta_{j}$ is called index of $J_{j}$, and we call invariant space $P_{j}$ neutral, of the first type, or of the second type, $n_{+}(J_{j})-n_{-}(J_{j})$ is equal to respectively zero, one or negative one.

Denote by $S_{+}$ ($S_{-}$) a matrix whose columns span the $R$-invariant subspace associated to eigenvalues with positive (respectively, negative) real parts. Then $S_{+}$ and $S_{-}$ are $J$-orthogonal to $S_{j}$ for every $j=1,\ldots,m$. Besides, all columns of matrices $S_{-}$, $S_{+}$, $S_{1},\ldots,S_{m}$ present a basis of ${\bf C}^{2n}$.

Fix $j\in \{1,\ldots,m\}$ and consider a matrix $V$ such that $JV$ is orthogonal to $S_{-}$, $S_{+}$, and all matrices $S_{k}$ with $k\ne j$. Then 

$$ \begin{array}{c} \det(R + V(JV)^{*} -\lambda I)=\det(R-\lambda I) (1+(JV)^{*}(R-\lambda I)^{-1}V)= \\
= \det(R-\lambda I) (1+(JV)^{*}S_{j}(J_{j}-\lambda I)^{-1}(S_{j}^{*}JS_{j})^{-1}S_{j}^{*}JV)
\end{array}
$$

Denote $P_{j}=S_{j}^{*}JS_{j}/\epsilon_{j}$.
Consider the matrix of quadratic form in the last brackets. We have $P_{j}^{-1}=(-1)^{m_{j}-1}P_{j}$, and 

$$ (J_{j}-\lambda I)^{-1} = 
\left(\begin{array}{cccc}
\frac{1}{i\omega_{j}-\lambda} & \frac{-1}{(i\omega_{j}-\lambda)^{2}} & \ldots & \frac{(-1)^{m_{j}-1}}{(i\omega_{j}-\lambda)^{m_{j}}} \\
0 & \frac{1}{i\omega_{j}-\lambda} & \ldots & \frac{(-1)^{m_{j}-2}}{(i\omega_{j}-\lambda)^{m_{j}-1}} \\
\ldots & \ldots & \ldots & \ldots \\ 0 & 0 & \ldots & 
\frac{1}{i\omega_{j}-\lambda} \end{array}\right).$$

Therefore 

$$ (J_{j}-\lambda I)^{-1}(\epsilon_{j}P_{j})^{-1} = 
\epsilon_{j}^{-1}\left(\begin{array}{cccc}
\frac{(-1)^{m_{j}}}{(i\omega_{j}-\lambda)^{m_{j}}} & 
\frac{(-1)^{m_{j}}}{(i\omega_{j}-\lambda)^{m_{j}-1}} & \ldots & \frac{(-1)^{m_{j}}}{i\omega_{j}-\lambda} \\ 
\frac{(-1)^{m_{j}-1}}{(i\omega_{j}-\lambda)^{m_{j}-1}} & 
\frac{(-1)^{m_{j}-1}}{(i\omega_{j}-\lambda)^{m_{j}-2}} & \ldots & 0 \\
\ldots & \ldots & \ldots & \ldots  \\
\frac{-1}{i\omega_{j}-\lambda} & 0 & \ldots & 0 \end{array}\right).
$$

For even $m_{j}=2r_{j}$ we have $\epsilon_{j}=(-1)^{r_{j}}\beta_{j}$, and 

$$ (J_{j}-\lambda I)^{-1}(\epsilon_{j}P_{j})^{-1} = 
\beta (-1)^{r_{j}} \left(\begin{array}{cccc}
\frac{1}{(i\omega_{j}-\lambda)^{2r_{j}}} & 
\frac{1}{(i\omega_{j}-\lambda)^{2r_{j}-1}} & \ldots & \frac{1}{i\omega_{j}-\lambda} \\ 
\frac{-1}{(i\omega_{j}-\lambda)^{2r_{j}-1}} & 
\frac{-1}{(i\omega_{j}-\lambda)^{2r_{j}-2}} & \ldots & 0 \\
\ldots & \ldots & \ldots & \ldots  \\
\frac{-1}{i\omega_{j}-\lambda} & 0 & \ldots & 0 \end{array}\right).
$$

Therefore for $\lambda = i\omega$ we have

$$ (J_{j}-i\omega I)^{-1}(\epsilon_{j}P_{j})^{-1} = 
\beta \left(\begin{array}{cccc}
\frac{1}{(\omega_{j}-\omega)^{2r_{j}}} & 
\frac{i}{(\omega_{j}-\omega)^{2r_{j}-1}} & \ldots & \frac{i(-1)^{r_{j}+1}}{\omega_{j}-\omega} \\ 
\frac{-i}{(\omega_{j}-\omega)^{2r_{j}-1}} & 
\frac{1}{(\omega_{j}-\omega)^{2r_{j}-2}} & \ldots & 0 \\
\ldots & \ldots & \ldots & \ldots  \\
\frac{i(-1)^{r_{j}}}{\omega_{j}-\omega} & 0 & \ldots & 0 \end{array}\right).
$$

If we choose vector $V$ such that $S_{k}^{*}(JV)=0$ for all $k\ne j$, and $S_{j}^{*}(JV)=(\delta,0,\ldots,0)^{*}$ (that is, $v$ is parallel to the last column of matrix $S_{j}$), then 

$$ 1+(JV)^{*}S_{j}(J_{j}-\lambda I)^{-1}(S_{j}^{*}JS_{j})^{-1}S_{j}^{*}JV = 
1+\delta^{2}\beta_{j} \frac{1}{(\omega_{j}-\omega)^{2r_{j}}}.
$$

Hence, for sufficiently small $\delta$ and positive $\beta_{j}$ all pure imaginary eigenvalues of matrix $R$ corresponding to block $J_{j}$ leave imaginary axis and become pairs of eigenvalues with non zero real parts symmetric with respect to imaginary axis.

If $\beta_{j}=-1$ then for sufficiently small $\delta$ exactly two eigenvalues remain on the imaginary axis, and the rest $m_{j}-2$ eigenvalues leave this axis. Among two remaining eigenvalues the bigger is of the first type, and the smaller is of the second type. With $\delta$ increasing the first eigenvalue goes up the imaginary axis, and the second eigenvalue goes down.

Now consider the blocks with odd dimension: $m_{j}=2r_{j}+1$.  Then $\epsilon_{j}=(-1)^{r_{j}}i\beta$, and 

$$ (J_{j}-\lambda I)^{-1}(\epsilon_{j}P_{j})^{-1} = 
\beta (-1)^{r_{j}} \left(\begin{array}{cccc}
\frac{i}{(i\omega_{j}-\lambda)^{2r_{j}+1}} & 
\frac{i}{(i\omega_{j}-\lambda)^{2r_{j}}} & \ldots & \frac{i}{i\omega_{j}-\lambda} \\ 
\frac{-i}{(i\omega_{j}-\lambda)^{2r_{j}}} & 
\frac{-i}{(i\omega_{j}-\lambda)^{2r_{j}-1}} & \ldots & 0 \\
\ldots & \ldots & \ldots & \ldots  \\
\frac{i}{i\omega_{j}-\lambda} & 0 & \ldots & 0 \end{array}\right).
$$

For $\lambda=i\omega$ we get

$$ (J_{j}-i\omega I)^{-1}(\epsilon_{j}P_{j})^{-1} = 
\beta \left(\begin{array}{cccc}
\frac{1}{(\omega_{j}-\omega)^{2r_{j}+1}} & 
\frac{i}{(\omega_{j}-\omega)^{2r_{j}}} & \ldots & \frac{(-1)^{r_{j}}}{\omega_{j}-\omega} \\ 
\frac{-i}{(\omega_{j}-\omega)^{2r_{j}}} & 
\frac{1}{(\omega_{j}-\omega)^{2r_{j}-1}} & \ldots & 0 \\
\ldots & \ldots & \ldots & \ldots  \\
\frac{(-1)^{r_{j}}}{\omega_{j}-\omega} & 0 & \ldots & 0 \end{array}\right).
$$

Again, if $V$ is such that $S_{k}^{*}(JV)=0$ for all $k\ne j$, and $S_{j}^{*}(JV)=(\delta,0,\ldots,0)^{*}$ (that is, $v$ is parallel to the last column of matrix $S_{j}$), then 

$$ 1+(JV)^{*}S_{j}(J_{j}-\lambda I)^{-1}(S_{j}^{*}JS_{j})^{-1}S_{j}^{*}JV = 
1+\delta^{2}\beta_{j} \frac{1}{(\omega_{j}-\omega)^{2r_{j}+1}}.
$$

Therefore, for small $\delta$ all eigenvalues of matrix $R$ corresponding to block $J_{j}$ but one leave imaginary axis. If $\beta=1$, then the eigenvalue which remains on imaginary axis is of the first type, and it goes up the imaginary axis with increasing positive $\delta$. If $\beta=-1$, then this eigenvalue is of the second type and it goes down imaginary axis with increasing positive $\delta$. 

Thus, we get the following result.

\begin{theorem}
There exists a nonnegative matrix $M$ such that the eigenvalues of matrix $R-tMJ$ with number $t$ increasing from zero have the following behaviour.

(a) For each Jordan block $J_{j}$ of matrix $R$ of odd dimension $2r_{j}+1$ with index $\beta=1$ exactly $2r_{j}$ eigenvalues leave imaginary axis with increasing $t$ from zero, and the rest eigenvalue goes up imaginary axis, and has the first type.

(b) For each Jordan block $J_{j}$ of matrix $R$ of odd dimension $2r_{j}+1$ with index $\beta=-1$ exactly $2r_{j}$ eigenvalues leave imaginary axis with increasing $t$ from zero, and the rest eigenvalue goes down imaginary axis, and has the second type.

(c) For each Jordan block $J_{j}$ of matrix $R$ of even dimension $2r_{j}$ with index $\beta=1$ all eigenvalues leave imaginary axis with increasing $t$ from zero.

(d) For each Jordan block $J_{j}$ of matrix $R$ of even dimension $2r_{j}$ with index $\beta=-1$ exactly $2r_{j}-2$ eigenvalues leave imaginary axis with increasing $t$ from zero. One of the rest eigenvalues goes up imaginary axis and has the first type, and the other eigenvalue goes down imaginary axis and has the second type.
\end{theorem}

From the representation above we see, that if two pure imaginary eigenvalues of $R-tMJ$ having different types meet for some number $t$, then further we do not change them. More precisely, we modify matrix $M$ (by eliminating the terms corresponding to Jordan blocks with these eigenvalues) such that further these eigenvalues and corresponding eigenvectors remain constant. 

We can finally get Hamiltonian matrix with all pure imaginary eigenvalues having even multiplicity only if each eigenvalue of the first type eventually meets an eigenvalue of the second type. Taking into account the fact that with increasing $t$ eigenvalues of the first type go up imaginary axis, and eigenvalues of the first type go down imaginary axis, it can happen if and only if for every $\omega$ the number of Jordan blocks of $R$ with eigenvalues $i\omega_{j}< \omega$ of the first type is not less than the number of Jordan blocks of $R$ with eigenvalues $i\omega_{j}< \omega$ of the second type plus the number of Jordan blocks of $R$ of neutral type with eigenvalue $i\omega$ and index $\beta=-1$. 

Denote 

$$ s(\omega) = m_{+}(\omega)- m_{-}(\omega) - m_{0}(\omega),$$

\noindent where $m_{+}(\omega)$ is the number of odd dimensional Jordan blocks of matrix $R$ with eigenvalues $ i\omega_{j}$ such that $\omega_{j}<\omega$, $\beta=1$,  $m_{-}(\omega)$ is the number of odd dimensional Jordan blocks of matrix $R$ with eigenvalues $ i\omega_{j}$ such that $\omega_{j}\le \omega$, $\beta=-1$, and $m_{0}(\omega)$ is the number of even dimensional Jordan blocks of matrix $R$ with eigenvalues $ i\omega$ such that $\beta=-1$.

\begin{theorem}
There exists a positive definite matrix $M$ such that matrix $R-MJ$ has no pure imaginary eigenvalues if and only if for every pure imaginary eigenvalue $i\omega_{j}$ of matrix $R$ we have

\begin{equation}\label{main}  
s(\omega_{j})\ge 0.
\end{equation}

\end{theorem}

Taking into account the result of the previous section, we get the following main result.

\begin{theorem}
The Riccati inequality (\ref{1}) has a solution if and only if inequality (\ref{main}) holds for all $\omega\in R$.
\end{theorem}

%
%
\section{A Numerical Example}

Consider the Riccati inequality:

\begin{equation}\label{ex} 
HA + A^{*} H + G_{1} - H B\Gamma ^{-1}B^{*}H < 0,
\end{equation}

where:

\noindent \begin{equation}
\{A, B, G \} = \left\{\left(
\begin{array}{ccc}
 1 & -1 & 1 \\
 0 & 1 & 1 \\
 0 & 0 & 1 \\
\end{array}
\right), \left(
\begin{array}{cc}
 1 & 0 \\
 1 & 0 \\
 0 & 1 \\
\end{array}
\right), \left(
\begin{array}{ccc}
 6 & -2 & -2 \\
 -2 & -3 & -2 \\
 -2 & -2 & -3.9 \\
\end{array}
\right)\right\}
\end{equation}
 and 
\[
\Gamma =\left (\begin{array}{cc}
 -10 & 0 \\
 0 & 0.1 \\
\end{array}\right)
\]


Now computing the eigenvalues of $A$ yields \{1,1,1\}. This it checks that matrix $A$ has no purely imaginary eigenvalues. Notice also that matrix $G$ and $\Gamma$ are Hermitian, and matrix $\Gamma$ is sign indefinite.

\noindent Associated Hamiltonian matrix of the inequality being:

\begin{equation} 
R = \left(\begin{array}{cc}
A&-B\Gamma ^{-1}B^{*} \\
 -G &-A^{*}\end{array}\right) = 
 \left(\begin{array}{cccccc}
 1 & -1 & 1 & 0.1 & 0.1 & 0. \\
 0 & 1 & 1 & 0.1 & 0.1 & 0. \\
 0 & 0 & 1 & 0. & 0. & -10. \\
 -6 & 2 & 2 & -1 & 0 & 0 \\
 2 & 3 & 2 & 1 & -1 & 0 \\
 2 & 2 & 3.9 & -1 & -1 & -1 \\
\end{array}\right).
\end{equation}

\noindent Notice that eigenvalues of $R$ are $\pm 6.0506i , \pm 1.5866i$ and $\pm 1.7964$. Thus, there are simple pure imaginary eigenvalues of matrix $R$.

The corresponding Riccati equation has no solution.

 \noindent Now consider some positive definite matrix, for instance 
 
\begin{equation}
\zeta =\left(
\begin{array}{ccc}
  4 & 2 & 2 \\
 2 & 4 & 2 \\
 2 & 2 & 4 \\
\end{array}
\right)
\end{equation}  

\noindent Then 

\begin{equation} 
G_1 = G + \zeta = 
\left( 
\begin{array}{ccc}
 10 & 0 & 0 \\
 0 & 1 & 0 \\
 0 & 0 & 0.1 \\
\end{array}
\right)
\end{equation}  

\noindent Notice that the matrix $\zeta$ is positive definite because it has the following eigenvalues 
$\{ 8, 2, 2 \}$. 

The equality becomes:

\begin{equation}\label{ex} 
HA + A^{*} H + G_{1} - H B\Gamma ^{-1}B^{*}H = 0
\end{equation}

With all conditions satisfied, a Hermitian matrix $H$ that satisfies (26) is equal to

\begin{equation}\label{H}
H=\left(
\begin{array}{ccc}
 128.485 & -178.389 & -7.18338 \\
 -178.389 & 259.987 & 12.4241 \\
 -7.18338 & 12.4241 & 1.25879 \\
\end{array}
\right)
\end{equation}

It is however straight-forward to check that this matrix $H$ is a solution of inequality (21).

%
%
\section{Conclusion}

We consider the problem of solvability of Riccati inequalities with arbitrary nondegenerous quadratic form. The necessary and sufficient condition for solvability of these inequalities were obtained and expressed in terms of the associated Hamiltonian matrices. The proof of this result is based on special representations of Hamiltonian matrices. An extension of this result to matrix pencils is a subject of possible future investigations.

\bibliographystyle{plain}

\end{document}